\documentclass{amsart}

\usepackage{amsfonts,amssymb,verbatim,amsmath,amsthm,latexsym,textcomp,amscd}
\usepackage{latexsym,amsfonts,amssymb,epsfig,verbatim}
\usepackage{amsmath,amsthm,amssymb,latexsym,graphics,textcomp}
\usepackage{graphicx}
\usepackage{color}
\usepackage{url}

\input xy
\xyoption{all}

\begin{document}

\newtheorem{theorem}{Theorem}[section]
\newtheorem{prop}[theorem]{Proposition}
\newtheorem{lemma}[theorem]{Lemma}
\newtheorem{cor}[theorem]{Corollary}
\newtheorem{definition}[theorem]{Definition}
\newtheorem{conj}[theorem]{Conjecture}
\newtheorem{rmk}[theorem]{Remark}
\newtheorem{claim}[theorem]{Claim}
\newtheorem{defth}[theorem]{Definition-Theorem}

\newcommand{\boundary}{\partial}
\newcommand{\C}{{\mathbb C}}
\newcommand{\integers}{{\mathbb Z}}
\newcommand{\natls}{{\mathbb N}}
\newcommand{\ratls}{{\mathbb Q}}
\newcommand{\reals}{{\mathbb R}}
\newcommand{\proj}{{\mathbb P}}
\newcommand{\lhp}{{\mathbb L}}
\newcommand{\tube}{{\mathbb T}}
\newcommand{\cusp}{{\mathbb P}}
\newcommand\AAA{{\mathcal A}}
\newcommand\BB{{\mathcal B}}
\newcommand\CC{{\mathcal C}}
\newcommand\DD{{\mathcal D}}
\newcommand\EE{{\mathcal E}}
\newcommand\FF{{\mathcal F}}
\newcommand\GG{{\mathcal G}}
\newcommand\HH{{\mathcal H}}
\newcommand\II{{\mathcal I}}
\newcommand\JJ{{\mathcal J}}
\newcommand\KK{{\mathcal K}}
\newcommand\LL{{\mathcal L}}
\newcommand\MM{{\mathcal M}}
\newcommand\NN{{\mathcal N}}
\newcommand\OO{{\mathcal O}}
\newcommand\PP{{\mathcal P}}
\newcommand\QQ{{\mathcal Q}}
\newcommand\RR{{\mathcal R}}
\newcommand\SSS{{\mathcal S}}
\newcommand\TT{{\mathcal T}}
\newcommand\UU{{\mathcal U}}
\newcommand\VV{{\mathcal V}}
\newcommand\WW{{\mathcal W}}
\newcommand\XX{{\mathcal X}}
\newcommand\YY{{\mathcal Y}}
\newcommand\ZZ{{\mathcal Z}}
\newcommand\CH{{\CC\HH}}
\newcommand\PEY{{\PP\EE\YY}}
\newcommand\MF{{\MM\FF}}
\newcommand\RCT{{{\mathcal R}_{CT}}}
\newcommand\PMF{{\PP\kern-2pt\MM\FF}}
\newcommand\FL{{\FF\LL}}
\newcommand\PML{{\PP\kern-2pt\MM\LL}}
\newcommand\GL{{\GG\LL}}
\newcommand\Pol{{\mathcal P}}
\newcommand\half{{\textstyle{\frac12}}}
\newcommand\Half{{\frac12}}
\newcommand\Mod{\operatorname{Mod}}
\newcommand\Area{\operatorname{Area}}
\newcommand\ep{\epsilon}
\newcommand\hhat{\widehat}
\newcommand\Proj{{\mathbf P}}
\newcommand\U{{\mathbf U}}
 \newcommand\Hyp{{\mathbf H}}
\newcommand\D{{\mathbf D}}
\newcommand\Z{{\mathbb Z}}
\newcommand\R{{\mathbb R}}
\newcommand\Q{{\mathbb Q}}
\newcommand\E{{\mathbb E}}
\newcommand\til{\widetilde}
\newcommand\length{\operatorname{length}}
\newcommand\tr{\operatorname{tr}}
\newcommand\gesim{\succ}
\newcommand\lesim{\prec}
\newcommand\simle{\lesim}
\newcommand\simge{\gesim}
\newcommand{\simmult}{\asymp}
\newcommand{\simadd}{\mathrel{\overset{\text{\tiny $+$}}{\sim}}}
\newcommand{\ssm}{\setminus}
\newcommand{\diam}{\operatorname{diam}}
\newcommand{\pair}[1]{\langle #1\rangle}
\newcommand{\T}{{\mathbf T}}
\newcommand{\inj}{\operatorname{inj}}
\newcommand{\pleat}{\operatorname{\mathbf{pleat}}}
\newcommand{\short}{\operatorname{\mathbf{short}}}
\newcommand{\vertices}{\operatorname{vert}}
\newcommand{\collar}{\operatorname{\mathbf{collar}}}
\newcommand{\bcollar}{\operatorname{\overline{\mathbf{collar}}}}
\newcommand{\I}{{\mathbf I}}
\newcommand{\tprec}{\prec_t}
\newcommand{\fprec}{\prec_f}
\newcommand{\bprec}{\prec_b}
\newcommand{\pprec}{\prec_p}
\newcommand{\ppreceq}{\preceq_p}
\newcommand{\sprec}{\prec_s}
\newcommand{\cpreceq}{\preceq_c}
\newcommand{\cprec}{\prec_c}
\newcommand{\topprec}{\prec_{\rm top}}
\newcommand{\Topprec}{\prec_{\rm TOP}}
\newcommand{\fsub}{\mathrel{\scriptstyle\searrow}}
\newcommand{\bsub}{\mathrel{\scriptstyle\swarrow}}
\newcommand{\fsubd}{\mathrel{{\scriptstyle\searrow}\kern-1ex^d\kern0.5ex}}
\newcommand{\bsubd}{\mathrel{{\scriptstyle\swarrow}\kern-1.6ex^d\kern0.8ex}}
\newcommand{\fsubeq}{\mathrel{\raise-.7ex\hbox{$\overset{\searrow}{=}$}}}
\newcommand{\bsubeq}{\mathrel{\raise-.7ex\hbox{$\overset{\swarrow}{=}$}}}
\newcommand{\tw}{\operatorname{tw}}
\newcommand{\base}{\operatorname{base}}
\newcommand{\trans}{\operatorname{trans}}
\newcommand{\rest}{|_}
\newcommand{\bbar}{\overline}
\newcommand{\UML}{\operatorname{\UU\MM\LL}}
\newcommand{\EL}{\mathcal{EL}}
\newcommand{\tsum}{\sideset{}{'}\sum}
\newcommand{\tsh}[1]{\left\{\kern-.9ex\left\{#1\right\}\kern-.9ex\right\}}
\newcommand{\Tsh}[2]{\tsh{#2}_{#1}}
\newcommand{\qeq}{\mathrel{\approx}}
\newcommand{\Qeq}[1]{\mathrel{\approx_{#1}}}
\newcommand{\qle}{\lesssim}
\newcommand{\Qle}[1]{\mathrel{\lesssim_{#1}}}
\newcommand{\simp}{\operatorname{simp}}
\newcommand{\vsucc}{\operatorname{succ}}
\newcommand{\vpred}{\operatorname{pred}}
\newcommand\fhalf[1]{\overrightarrow {#1}}
\newcommand\bhalf[1]{\overleftarrow {#1}}
\newcommand\sleft{_{\text{left}}}
\newcommand\sright{_{\text{right}}}
\newcommand\sbtop{_{\text{top}}}
\newcommand\sbot{_{\text{bot}}}
\newcommand\sll{_{\mathbf l}}
\newcommand\srr{_{\mathbf r}}
\newcommand\geod{\operatorname{\mathbf g}}
\newcommand\mtorus[1]{\boundary U(#1)}
\newcommand\A{\mathbf A}
\newcommand\Aleft[1]{\A\sleft(#1)}
\newcommand\Aright[1]{\A\sright(#1)}
\newcommand\Atop[1]{\A\sbtop(#1)}
\newcommand\Abot[1]{\A\sbot(#1)}
\newcommand\boundvert{{\boundary_{||}}}
\newcommand\storus[1]{U(#1)}
\newcommand\Momega{\omega_M}
\newcommand\nomega{\omega_\nu}
\newcommand\twist{\operatorname{tw}}
\newcommand\modl{M_\nu}
\newcommand\MT{{\mathbb T}}
\newcommand\Teich{{\mathcal T}}
\renewcommand{\Re}{\operatorname{Re}}
\renewcommand{\Im}{\operatorname{Im}}

\title{Uniformization of higher genus finite type Log-riemann surfaces}

\author{Kingshook Biswas}
\address{RKM Vivekananda University, Belur Math, WB-711 202, India}

\author{Ricardo Perez-Marco}
\address{CNRS, LAGA, UMR 7539, Universit\'e
Paris 13, Villetaneuse, France}

\begin{abstract}
We consider a log-Riemann surface $\mathcal{S}$ with a finite
number of ramification points and finitely generated fundamental
group. The log-Riemann surface is equipped with a local
holomorphic difffeomorphism $\pi : \mathcal{S} \to \C$. We prove
that $\mathcal{S}$ is biholomorphic to a compact Riemann surface with
finitely many punctures $S$, and the pull-back of the $1$-form $d\pi$ under
the biholomorphic map $\phi : S \to \mathcal{S}$ is a $1$-form $\omega = \phi^* d\pi$
with isolated singularities at the punctures of exponential type, i.e. near each puncture
$p$, $\omega = e^h \cdot \omega_0$ where $h$ is a function meromorphic near $p$ and
$\omega_0$ a $1$-form meromorphic near $p$.
\end{abstract}

\bigskip

\maketitle

\tableofcontents

\section{Introduction}

\medskip

In \cite{bipm1} we defined the notion of log-Riemann surface, as a
Riemann surface $\SSS$ equipped with a local diffeomorphism $\pi :
\SSS \to \C$ such that the set of points $\RR$ added in the
completion ${\SSS}^* = \SSS \sqcup \RR$ of $\SSS$ with respect to
the flat metric on $\SSS$ induced by $|d\pi|$ is discrete. It was shown in
\cite{bipm1} that $\pi$ extends to the points $p \in \RR$, and is a covering
of a punctured neighbourhood of $p$ onto a punctured disk in $\C$;
the point $p$ is called a ramification point of $\SSS$ of order
equal to the degree of the covering $\pi$ near $p$. The finite
order ramification points may be added to $\SSS$ to give a Riemann
surface $\SSS^\times$, called the finite completion of $\SSS$. We call a
log-Riemann surfaces of {\it finite type} if it has finitely many ramification
points and finitely generated fundamental group.

\medskip

In \cite{bipm2} we considered
log-Riemann surfaces of finite type with simply connected finite completion
and showed that these surfaces were biholomorphic to $\C$, with uniformization given by an entire
function of the form $F(z) = \int Q(z)e^{P(z)} dz$ for some polynomials $P, Q$, in the
sense the entire function $F = \pi \circ \phi$, where  $\phi : \C \to \SSS^{\times}$ is
a biholomorphic map, is of this form. It was also shown that conversely
any such entire function arises as the uniformization of a simply
connected finite type log-Riemann surface. The entire functions of the above form
have also been studied by Nevanlinna \cite{nevanlinna1} and
M. Taniguchi \cite{taniguchi1}.

\medskip

The aim of this article is to generalize these uniformization
theorems to arbitrary finite type log-Riemann surfaces. The role of the entire
functions above is played by functions $f$ with finitely many critical points
on punctured surfaces such that $df$ has {\it exponential singularities} at the
punctures; here we say an isolated singularity $p$
of a holomorphic $1$-form $\omega$ is an exponential singularity if locally
$\omega = e^{h} \cdot \omega_0$ where $h$ is a function with a pole at $p$ and
$\omega_0$ a $1$-form meromorphic near $p$.

\medskip

Ramification points of log-Riemann surfaces are defined in the
next section; we define:

\medskip

\noindent {\bf Definition.} A log-Riemann surface $\SSS$ is of
finite type if it has finitely many ramification points and
finitely generated fundamental group.

\medskip

We prove the following:

\medskip

\begin{theorem} \label{uniformthm} Let $\SSS$ be a log-Riemann
surface of finite type. Then there is a compact Riemann surface $C$, a finite
set $A \subset C$ and a biholomorphic map
$\phi : C - A \to \SSS^{\times}$ such that the map
$f = \pi \circ \phi$ has finitely many critical points, and
$df$ has exponential singularities or poles at the points of $A$.
\end{theorem}

\medskip

Conversely we have:

\begin{theorem} \label{conversethm} Let $C$ be a compact Riemann
surface, $A \subset C$ a finite set and $f$ a holomorphic
function on $C - A$ with finitely many critical points such that $df$
has either exponential singularities or poles at the points of $A$. Then there
is a finite type log-Riemann surface $\SSS$ and a biholomorphic
map $\phi : C - A \to \SSS^{\times}$ such that $f = \pi \circ \phi$.
\end{theorem}

\medskip

The proof of Theorem \ref{uniformthm} proceeds in outline as
follows: we first show that $\SSS^{\times}$ has finitely many ends,
each of which is homeomorphic to a punctured disk.
We then show that each end can be isometrically embedded in a simply connected finite type
log-Riemann surface. It follows from the uniformization theorem for simply
connected log-Riemann surfaces of \cite{bipm2} that each end is biholomorphic to a
punctured disk, and moreover the $1$-form $d\pi$ has an
exponential singularity at the puncture. Since $\SSS^{\times}$ has
finitely many ends each of which is biholomorphic to a punctured
disk, it follows that $\SSS^{\times}$ is biholomorphic to a
compact Riemann surface minus finitely many punctures, and $d\pi$ has exponential
singularities at the punctures.

\medskip

For the converse Theorem \ref{conversethm}, since $f$ is a local diffeomorphism
away from the finite set of critical points it suffices to show
that there are only finitely many points added when $C$ is completed with
respect to the flat path metric induced by $|df|$. Using the fact
that $df$ has exponential singularities at the punctures, we show
that a neighbourhood of each puncture is bi-Lipschitz equivalent
to a neighbourhood of infinity in $\C$ under the path metric
induced by $|Q(z) e^{P(z)} dz|$ for some polynomials $P, Q$. It
follows from the converse of the uniformization theorem of
\cite{bipm2} that only finitely many points are added in the completion of
a neighbourhood of each puncture of $C$, and hence in the completion of $C$.

\medskip

In section 2 we study the topology of finite type log-Riemann
surfaces, using the tools introduced in \cite{bipm2}, and show
that there are finitely many ends, which are all punctured disks.
Then in section 3 we prove Theorems \ref{uniformthm},
\ref{conversethm} above.

\bigskip

\section{Topology of finite type log-Riemann surfaces}

\medskip

We first recall some basic properties of log-Riemann surfaces from
\cite{bipm1}. A log-Riemann surface $\SSS$ is
equipped with a local holomorphic diffeomorphism $\pi : \SSS \to \C$ such that
the following holds: the flat metric $|d\pi|$ induces a path metric $d$ on
$\SSS$; letting $\overline{\SSS} = \SSS \sqcup \RR$ be the metric completion of $\SSS$,
the set $\RR$ of points added is discrete.

\medskip

Moreover for each $w^* \in \RR$ the map $\pi$
restricted to a sufficiently small punctured disk $B(w^*,
\epsilon) - \{w^*\}$ is a connected covering of a punctured disk
$B(\pi(w^*), \epsilon) - \{\pi(w^*)\}$ in $\C$. The point $w^*$ is
called a ramification point of $\SSS$ of order equal to the degree
$1 \leq n \leq +\infty$ of the covering.

\medskip

A punctured disk neighbourhood of a finite order ramification
point is biholomorphic to a punctured disk, hence the finite
ramification points may be added to $\SSS$ to obtain a Riemann
surface $\SSS^{\times}$ called the {\it finite completion} of
$\SSS$.

\medskip

\subsection{Decomposition into stars}

\medskip

The map $\pi$ restricted to $\SSS$ is a local isometry. The locally Euclidean
metric on $\SSS$ can be used to define a decomposition into simply
connected open sets called {\it stars} defined as follows:

\medskip

Let $w_0 \in \SSS$. Given an angle $\theta \in
\reals/2\pi\integers$, there is a $0 < \rho(w_0,\theta) \leq
+\infty$ and a unique maximal unbroken geodesic
segment $\gamma(w_0, \theta) : [0, \rho(w_0, \theta)) \to \SSS$ starting
at $w_0$ which is the lift of the line segment $\{ \pi(w_0)
+ t e^{i\theta} : 0 \leq t < \rho(w_0,\theta) \}$, such that
$\gamma(w_0, \theta)(t) \to w^* \in \RR$ as $t \to \rho(w_0, \theta)$
if $\rho(w_0, \theta) < +\infty$.

\medskip

\begin{definition} The star of $w_0 \in \SSS$ is the union of all
maximal unbroken geodesics starting at $w_0$,
$$
V(w_0) := \bigcup_{\theta \in \reals/2\pi\integers} \gamma(w_0, \theta)
$$
Similarly we also define for a ramification point $w^*$ of order $n \leq +\infty$
the star $V(w^*)$ as the union of all maximal unbroken
geodesics $\gamma(w^*, \theta)$ starting from $w^*$, where the angle $\theta \in [-n\pi,
n\pi)$:
$$
V(w^*) := \{ \gamma(w^*, \theta)(t) : 0 \leq t < \rho(w^*, \theta) , -n\pi \leq \theta \leq n\pi \}
$$
\end{definition}

\medskip

The following Proposition is proved in \cite{bipm2}:

\begin{prop} \label{stars} For any $w_0 \in \SSS$ the star $V(w_0)$ is a simply connected open subset
of $\SSS$, and $\pi$ maps $V(w_0)$ biholomorphically onto a slit plane $\C - L$, where $L$ is a
locally finite union of closed half-lines.
The boundary $\partial V(w_0) \subset \SSS$ is a disjoint union of maximal unbroken geodesic
segments $\gamma : (0, +\infty) \to \SSS$ such that
$\gamma(t) \to w^* \in \RR$ as $t \to 0$, $\gamma(t) \to \infty$
as $t \to \infty$..
\end{prop}

\medskip

The set of ramification points $\RR$ is discrete, hence countable.
Let $L \supset \pi(\RR)$ be the union in $\C$ of all straight
lines joining points of $\pi(\RR)$. Then $\C - L$ is dense in $\C$.
By a {\it generic fiber} we mean a
fiber $\pi^{-1}(z_0) = \{ w_i \}$ of $\pi$ such that $z_0 \in \C - L$.

\medskip

From \cite{bipm2} we have:

\medskip

\begin{prop} \label{starsdense} Let $\{w_i\}$ be a generic fiber. Then:

\noindent{(1)} The stars $\{ V(w_i) \}$ are disjoint.

\noindent{(2)} For $w_i \neq w_j$, the components of $\partial V(w_i), \partial V(w_j)$ are
either disjoint or equal, and each component can belong to at most two such
stars.

\noindent{(3)} The union of the stars is dense in $\SSS$:
$$
\SSS = \overline{\bigcup_i V(w_i)} = \bigcup_i \overline{V(w_i)}
$$

\end{prop}

\medskip

The above Proposition gives a cell
decomposition of $\SSS$ into cells $V(w_i)$ glued along boundary
arcs $\gamma \subset \partial V(w_i) \cap \partial V(w_j)$.

\bigskip

\subsection{The skeleton and fundamental group}

\medskip

Let $\pi^{-1}(z_0) = \{w_i\}$ be a generic fiber. The $1$-skeleton of the
cell decomposition into stars gives an associated graph:

\medskip

\begin{definition} The {\it skeleton} $\Gamma(\SSS, z_0) = (\VV, \EE)$ is the graph with
vertices given by the stars $V_i = V(w_i)$, and an edge between $V(w_i)$
and $V(w_j)$ for each connected component $\gamma$ of
$\partial V(w_i) \cap \partial V(w_j)$. Each edge corresponds to a geodesic ray
$\gamma : (0, +\infty) \to \SSS$ starting at a ramification point.
This gives us a map from edges to ramification points,
$\hbox{foot} : \gamma \mapsto \hbox{foot}(\gamma)
:= \lim_{t \to 0} \gamma(t) \in \RR \cap \overline{V(w_i)} \cap
\overline{V(w_j)}$.

\medskip

To each ramification point $w^* \in \RR$ we associate the subgraph $\Gamma(w^*)$ with
vertices and edges
$$
\VV(w^*) := \{ V_i : w^* \in \overline{V_i} \},
\EE(w^*) := \{ \gamma : \hbox{foot}(\gamma) = w^* \}.
$$

\medskip

We note that if $\RR \neq \emptyset$ then
$$
\Gamma(\SSS, z_0) = \cup_{w^* \in \RR} \Gamma(w^*)
$$
\end{definition}

\medskip

From \cite{bipm2} we have the following Propositions:

\medskip

\begin{prop} If $w^*$ is of finite order $n$ then
$\Gamma(w^*)$ is a cycle in $\Gamma(\SSS, z_0)$ of
length $n$. If $w^*$ is of infinite order then $\Gamma(w^*)$
is a bi-infinite line in $\Gamma(\SSS, z_0)$.
\end{prop}

\medskip

\begin{prop} \label{retract1} The log-Riemann surface $\SSS$ deformation retracts onto $\Gamma(\SSS,z_0)$.
In particular $\pi_1(\SSS) = \pi_1(\Gamma(\SSS,z_0))$ is a free group.
\end{prop}

\medskip

The relation of $\Gamma(\SSS, z_0)$ to the finitely completed
log-Riemann surface $\SSS^{\times}$ is as follows:

\medskip

\begin{definition} The finitely completed skeleton
$\Gamma^{\times}(\SSS, z_0)$ is the graph obtained from
$\Gamma(\SSS, z_0)$ as follows: for each finite order ramification point
$w^*$, add a vertex $v = v(w^*)$ to $\Gamma(\SSS, z_0)$, remove all edges in
the cycle $\Gamma(w^*)$ and add an edge from $v_i$ to $v$ for each vertex $v_i$ in the
cycle $\Gamma(w^*)$.
\end{definition}

\medskip

We have:

\medskip

\begin{prop} \label{retract2} The finitely completed log-Riemann surface
$\SSS^{\times}$ deformation retracts onto the finitely completed
skeleton $\Gamma^{\times}(\SSS, z_0)$. In particular
$\pi_1(\SSS^{\times}) = \pi_1(\Gamma^{\times}(\SSS,z_0))$ is a
free group.
\end{prop}

\medskip

Finally we have:

\medskip

\begin{prop} \label{fingen} If $\pi_1(\SSS)$ is finitely generated then $\Gamma(\SSS,
z_0)$ is of the form
$\Gamma_0 \cup \TT$ where $\Gamma_0$ is a finite connected subgraph and
$\TT = \sqcup_{i \in I} T_i$ is a finite disjoint
union of rooted trees $T_i$ each intersecting $\Gamma_0$ only at its
root. A similar decomposition $\Gamma^{\times}(\SSS, z_0) = \Gamma^{\times}_0 \cup {\TT}^{\times}$
holds if $\pi_1(\SSS^{\times})$ is finitely generated.
\end{prop}

\medskip

\noindent{\bf Proof:} Let $\Gamma_1 \subset \Gamma(\SSS, z_0)$ be
a maximal subtree. Then $\pi_1(\Gamma(\SSS,z_0))$ is a free group
with generators corresponding to edges of $\Gamma(\SSS, z_0) -
\Gamma_1$. It follows that $\Gamma(\SSS, z_0)$ is given by
adjoining finitely many edges $e_1, \dots, e_n$ to a tree $\Gamma_1$.
Let $\Gamma'_1 \subset \Gamma_1$ be a finite subtree containing the
vertices of $e_1, \dots, e_n$ and $\Gamma_0$ the finite connected subgraph
given by adjoining $e_1, \dots, e_n$ to $\Gamma'_1$. Then the closure $T$
of each connected component of $\Gamma(\SSS, z_0) - \Gamma_0$ is
contained in $\Gamma_1$, hence is a tree. Moreover since $\Gamma'_1$ is connected,
$T$ only intersects $\Gamma_0$ at a single vertex $v$ (otherwise $T \cup \Gamma'_1 \subset \Gamma_1$
would contain a notrivial loop). We let $T(v)$ be the union of all such
trees $T'$ intersecting $\Gamma_0$ at $v$. Then $T(v)$ is a tree intersecting $\Gamma_0$
only at $v$, the trees $T(v), v $ a vertex of $\Gamma_0$, are disjoint and there are finitely
many such. Letting $\TT = \sqcup_{v \in \Gamma_0} T(v)$, we have
$\Gamma(\SSS, z_0) = \Gamma_0 \cup \TT$ as required. The proof for $\SSS^{\times},
\Gamma^{\times}(\SSS, z_0)$ is similar. $\diamond$

\medskip

\subsection{Ends of finite type log-Riemann surfaces}

%
%
%

\medskip

Let $\SSS$ be a finite type log-Riemann surface, i.e. $\SSS$ has
finitely many ramification points and finitely generated
fundamental group. We fix a generic fiber $\pi^{-1}(z_0)$ and
the associated skeleton $\Gamma = \Gamma(\SSS, z_0)$.
Let $\RR_{\infty} \subset \RR$ be the set of infinite ramification points of $\SSS$.

\medskip

\begin{prop} \label{halflines} The skeleton $\Gamma$ is of the form
$$
\Gamma = \Gamma_0 \sqcup \left(\sqcup_{w^* \in \RR_{\infty}} L^+(w^*) \sqcup L^-(w^-)\right)
$$
where $\Gamma_0$ is a finite connected subgraph
such that $(\Gamma(w^*) - \Gamma_0)$ is a disjoint union of two
open half-lines $L^+(w^*), L^-(w^*)$ for all $w^* \in \RR_{\infty}$.
\end{prop}

\medskip

\noindent{\bf Proof:} By Proposition
\ref{fingen}, there is a finite connected subgraph $\Gamma_0$ such that
$\Gamma = \Gamma_0 \cup \TT$ where $\TT$ is a finite
disjoint union of trees $\TT = \sqcup_{i \in I}
T_i$ each meeting $\Gamma_0$ in a single vertex. In particular,
all cycles of $\Gamma$ are contained in $\Gamma_0$.
Hence the subgraphs $\Gamma(w^*)$ associated to the
finite ramification points $w^*$ of $\SSS$ are contained in
$\Gamma_0$, and recalling that $\Gamma = \cup_{w^* \in \RR}
\Gamma(w^*)$, it follows that $\Gamma = \Gamma_0 \cup
(\cup_{w^* \in \RR_{\infty}} \Gamma(w^*))$.

\medskip

Each intersection
$\Gamma(w^*) \cap \Gamma_0, w^* \in \RR_{\infty}$, if nonempty, is connected, since
otherwise there would be a finite subinterval of $\Gamma(w^*)$ contained in
$\TT$ meeting $\Gamma_0$ in exactly two vertices, which is not possible. Hence
either $\Gamma(w^*)$ is disjoint from $\Gamma_0$ or
$\Gamma(w^*) \cap \Gamma_0$ is connected and $(\Gamma(w^*) - \Gamma_0)$
is a disjoint union of two open half-lines $L^+(w^*), L^-(w^*)$.

\medskip

If $\Gamma(w^*)$ is disjoint from $\Gamma_0$
then $\Gamma(w^*)$ is contained in a tree $T$ meeting
$\Gamma_0$ in a single vertex $v$ say,
and there is a simple arc $\gamma$ in
$T$ starting at $v$ which meets $\Gamma(w^*)$ in a single
vertex. It follows that by adding such paths $\gamma$ to
$\Gamma_0$, we can ensure that $\Gamma_0$ meets each
$\Gamma(w^*), w^* \in \RR_{\infty}$ and each $\Gamma(w^*) - \Gamma_0$ is a disjoint union
of two open half-lines $L^+(w^*), L^-(w^*)$ as required.

\medskip

It remains to check that for $w^*_1 \neq w^*_2$, any two half-lines
$L^{\pm}(w^*_1), L^{\pm}(w^*_2)$ are disjoint; this follows from observing that any
intersection would give rise to a cycle of $\Gamma$ not contained in $\Gamma_0$. $\diamond$

\medskip

We make the convention of labelling the half-lines $L^+(w^*), L^-(w^*)$
above such that any branch of $\arg(w - \pi(w^*))$ defined in a punctured
neighbourhood of $w^*$ is bounded below in
the stars of $\SSS$ corresponding to $L^+(w^*)$ and is bounded above
in those corresponding to $L^-(w^*)$. We let $l(w^*) \subset \C$ be
the closed half-line in the direction $(\pi(w^*) - z_0)$ starting
from $\pi(w^*)$. An immediate corollary of the
above Proposition is the following:

\medskip

\begin{lemma} \label{cleansheets} ('Clean sheets') For any star $V
\in L^{\pm}(w^*), w^* \in \RR_{\infty}$,
we have $\overline{V} \cap \RR = \{w^*\}$ (closure taken in $\overline{\SSS}$)
and $\pi$ maps $V$ biholomorphically onto the slit plane $\C - l(w^*)$.
\end{lemma}

\medskip

\noindent{\bf Proof:} The half-lines $L^{\pm}(w^*), L^{\pm}(w^*_1)$ are disjoint for
$w^*_1 \neq w^*$, so $V \notin \Gamma(w^*_1)$ for
$w^*_1 \neq w^*$ and $\overline{V} \cap \RR = \{w^*\}$, from which the second
conclusion follows. $\diamond$

\medskip

We now describe a neighbourhood of infinity in $\SSS^{\times}$:

\medskip

Fix $\epsilon > 0$ small enough so that the punctured
neighbourhoods $U(w^*) := \{ 0 < d(w, w^*_j) \leq \epsilon \}, w^* \in \RR_{\infty}$
are disjoint, and $\pi$ restricted to each $U(w^*)$ is a universal covering
of $\pi(U(w^*))$. Fix $R > |z_0|, |\pi(w^*)|+\epsilon, w^* \in \RR_{\infty}$, and
let $U(\infty) := \pi^{-1}(\{|z| \geq R
\})$. Then on each component of $U(\infty)$, $\pi$ is a connected covering of $\pi(U(\infty))$.
For $w^* \in \RR_{\infty}$, let $S^+(w^*)$ (respectively $S^-(w^*)$) be the
closure (in $\SSS^{\times}$) of the union of the stars in $L^+(w^*)$
(respectively $L^-(w^*)$). Let $S_0$ be the closure of the finite union of
the stars in $\Gamma_0$. Then by Proposition \ref{halflines}, we
have
$$
\SSS^{\times} = S_0 \cup \left(\cup_{w^* \in \RR_{\infty}} S^+(w^*) \cup
S^-(w^*)\right)
$$
Moreover, since $L^+(w^*), L^-(w^*)$ and $\Gamma_0$ are
connected, the sets $S^+(w^*), S^-(w^*)$ and $S_0$ are
connected.

\medskip

\begin{prop} \label{nbhdinfinity} The set
$$
E = U(\infty) \cup \left(\cup_{w^* \in \RR_{\infty}} U(w^*) \cup S^+(w^*) \cup
S^-(w^*)\right)
$$
is a neighbourhood of infinity in $\SSS^{\times}$, i.e.
$\SSS^{\times} - E$ is pre-compact.
\end{prop}

\medskip

\noindent{\bf Proof:} Let $(z_k)$ be a sequence in $\SSS^{\times} -
E$. Now for any star $V$ it is easy to see from Propositions \ref{stars}, \ref{starsdense}
that the set $V - (U(\infty) \cup_{w^* \in \RR_{\infty}} U(w^*))$ is precompact in
$\SSS^{\times}$. As $S_0$ is the closure of a finite union of stars,
$S_0 - (U(\infty) \cup_{w^* \in \RR_{\infty}} U(w^*))$ is precompact.
Since $z_k$ does not lie in any $S^{\pm}(w^*)$, $(z_k)$ must be
contained in $S_0$, and $z_k \notin (U(\infty) \cup_{w^* \in \RR_{\infty}} U(w^*))$,
hence the preceding remark implies that $(z_k)$ has a convergent subsequence. $\diamond$

\medskip

\begin{lemma} \label{ends1} The set $E$ has finitely many connected
components.
\end{lemma}

\medskip

\noindent{\bf Proof:} Since $S^+(w^*), S^-(w^*), S_0$ are
connected unions of closures of stars,
it is easy to see that the intersections $U(\infty) \cap S^+(w^*), U(\infty) \cap S^-(w^*),
U(\infty) \cap S_0$ are connected. As these intersections cover $U(\infty)$, any
component of $U(\infty)$ must be a finite union of these sets,
therefore $U(\infty)$ has finitely many components.
Similarly, the sets $U(w^*), S^{\pm}(w^*), w^* \in \RR_{\infty}$ are connected,
so any component of $E$ is a finite union of these sets
and of components of $U(\infty)$. It follows that $E$ has finitely many
components. $\diamond$

\medskip

\begin{lemma} \label{ends2} Let $U$ be a connected component of
$U(\infty)$. Let $\gamma : \reals \to \C$ be the curve $\gamma(t) = Re^{it}$. Then either:

\smallskip

\noindent (1) $U$ does not meet $\cup_{w^* \in \RR_{\infty}} (S^+(w^*) \cup S^-(w^*))$
and $\pi : U \to \{ |z| \geq R \}$ is a finite sheeted covering.

\smallskip

\noindent{or:}

\smallskip

\noindent (2) $U$ meets $\cup_{w^* \in \RR_{\infty}} (S^+(w^*) \cup S^-(w^*))$
and $\pi : U \to \{ |z| \geq R \}$ is a
universal covering. Moreover there are unique
$w^*_+ = w^*_+(U), w^*_- = w^*_-(U) \in \RR_{\infty}$ such that for any lift
$\tilde{\gamma}$ of $\gamma$ to $U$, $\tilde{\gamma}(t) \to \infty$ through
$S^+(w^*_+)$ as $t \to +\infty$ and $\tilde{\gamma}(t) \to \infty$ through $S^-(w^*_-)$ as $t \to
-\infty$.

\end{lemma}

\medskip

\noindent{\bf Proof:} We know $\pi : U \to \{ |z| \geq R \}$ is a covering.
Let $\tilde{\gamma}$ be a lift of $\gamma$ to $U$. As $\{|z| \geq R \}$ is a punctured
disk, clearly $\pi$ restricted to $U$ is a finite sheeted covering if and only if
$\tilde{\gamma}$ is periodic.

\medskip

Suppose $\tilde{\gamma}$ does not meet
$\cup_{w^* \in \RR_{\infty}} (S^+(w^*) \cup S^-(w^*))$. Then $\tilde{\gamma}$ is
contained in $\{ |\pi(w)| = R \} \cap S_0$ which is
compact (since $S_0$ is the closure of a finite union of stars),
so it follows that $\tilde{\gamma}$ is periodic and $\pi$ is a
finite sheeted covering. This proves (1).

\medskip

Now suppose $\tilde{\gamma}$ meets some $S^{\pm}(w^*)$, say $S^+(w^*)$. Then it follows from
Lemma \ref{cleansheets} that $\tilde{\gamma}(t) \to \infty$ through
$S^+_(w^*)$ as $t \to +\infty$, so $\pi$ is a universal covering.
Therefore $\tilde{\gamma}(t) \to \infty$
as $t \to -\infty$ as well. Since $\{ |\pi(w)| = R \} \cap S_0$ is compact,
it follows that as $t \to -\infty$, $\tilde{\gamma}(t) \to \infty$  through
$\cup_{w^* \in \RR_{\infty}} (S^+(w^*) \cup S^-(w^*))$. As these sets
are disjoint, $\tilde{\gamma}(t) \to \infty$ as $t \to -\infty$ through one of
these sets, say $S^{\pm}(\tilde{w}^*)$. As any branch of $\arg(w)$
is bounded below on $S^+(\tilde{w}^*)$, we must have $\tilde{\gamma}(t) \to
\infty$ through $S^-(\tilde{w}^*)$ as $t \to -\infty$. Letting $w^*_+ = w^*,
w^*_- = \tilde{w}^*$, this proves (2).
$\diamond$

\medskip

It follows from the Lemma that each $S^+(w^*), S^-(w^*),$ meets a
unique component of $U(\infty)$ which is necessarily as in case (2)
above; denote these components by $U^+(\infty, w^*)$ and $U^-(\infty,
w^*)$ respectively.

\medskip

We define maps $u, d: \RR_{\infty} \to \RR_{\infty}$ (for 'up', 'down'
respectively) by $u(w^*) := w^*_+(U^-(\infty, w^*)),
d(w^*) := w^*_-(U^+(\infty, w^*))$.
Since $w^* = w^*_-(U) = w^*_+(U')$, the maps $u,d$ are mutual
inverses, and so $\RR_{\infty}$ splits into disjoint cycles
$\RR_{\infty} = \RR_1 \sqcup \dots \sqcup \RR_m$
invariant under $u,d$.

\medskip

\begin{lemma} \label{cycles} For any cycle $\RR_i$, the set
$$
W(\RR_i) := \bigcup_{w^* \in \RR_i} \left(S^-(w^*) \cup U(w^*) \cup S^+(w^*) \cup U^-(\infty, w^*) \right)
$$
is connected. Moreover for distinct cycles $\RR_i, \RR_j$ the sets
$W(\RR_i), W(\RR_j)$ are disjoint.
\end{lemma}

\medskip

\noindent{\bf Proof:} Since, for any $w^*$, $S^-(w^*)$ and $S^+(u(w^*))$ meet $U^-(\infty,
w^*)$, it follows easily that $W(\RR_i)$ is connected.

\medskip

If for two cycles $\RR_i, \RR_j$ the sets $W(\RR_i), W(\RR_j)$
meet, then, for some $w^*_1 \in \RR_i, w^*_2 \in \RR_j$, the sets
$A_k = (S^-(w^*_k) \cup U(w^*_k) \cup S^+(w^*_k) \cup U^-(\infty, w^*_k)), k=1,2$,
must meet. Then either $w^*_1 = w^*_2$, or
$w^*_1 \neq w^*_2$, in which case one of the two intersections
$S^+(w^*_1) \cap U^-(\infty, w^*_2), S^+(w^*_2) \cap U^-(\infty,
w^*_1)$ must be nonempty (all other intersections between the constituents of
$A_1, A_2$ are empty), and consequently either $w^*_1 =
u(w^*_2)$ or $w^*_2 = u(w^*_1)$. Thus $\RR_i = \RR_j$. $\diamond$

\medskip

\begin{lemma} \label{ends3} Let $C$ be a component of $E$. Then either:

\smallskip

\noindent (1) $C = U$ for some component $U$ of $U(\infty)$ such that $\pi :
U \to \{ |z| \geq R \}$ is a finite-sheeted covering.

\smallskip

\noindent{or:}

\smallskip

\noindent (2) $C = W(\RR_i)$ for some cycle of infinite
ramification points $\RR_i \subset \RR_{\infty}$.
\end{lemma}

\medskip

\noindent{\bf Proof:} We first observe that any component $U$ of
$U(\infty)$ as in case (1) of Lemma \ref{ends2} does not meet
$\cup_{w^* \in \RR_{\infty}} (S^-(w^*) \cup U(w^*) \cup S^+(w^*))$,
and is hence a connected component of $E$, while any component of
$U(\infty)$ as in case (2) of Lemma \ref{ends2}
is contained in a unique $W(\RR_i)$ by Lemma \ref{cycles}.
We know $C$ is a finite union of components of
$U(\infty)$ and of the sets $(U(w^*) \cup S^+(w^*) \cup S^-(w^*))$.
In particular $C$ must meet $U(\infty)$.

\medskip

If $C$ contains a component $U$
of $U(\infty)$ as in case (1) of Lemma \ref{ends2} then by the previous
observation $C = U$ and we are done.

\medskip

Otherwise all components of $U(\infty)$
contained in $C$ are as in case (2) of Lemma \ref{ends2}, and are contained
in the union of the sets $W(\RR_i)$. As the sets $(U(w^*) \cup S^+(w^*) \cup S^-(w^*))$
are also contained in the union of the $W(\RR_i)$, it follows that $C$ is contained
in the union of the $W(\RR_i)$. Since the $W(\RR_i)$ are
pairwise disjoint, $C = W(\RR_i)$ for some $i$. $\diamond$

\medskip

\begin{lemma} \label{ends4} Each component of $E$ is
homeomorphic to a closed disk with a puncture, $\{ 0 < |z| \leq 1 \}$.
\end{lemma}

\medskip

\noindent{\bf Proof:} This is clear for the components $U$ of $U(\infty)$
which are finite sheeted coverings of $\{ |z| \geq R \}$. For a
component $W(\RR_i)$ with $\RR_i = \{w^*_0, \dots, w^*_{n-1} \}$,
where $w^*_{j+1} = u(w^*_j)$, we can decompose $W(\RR_i)$ as follows:

\medskip

Define $\tilde{U}(w^*_j) = \overline{U(w^*_j) - (S^-(w^*_j) \cup
S^+(w^*_j))}$ \hfill and  $\tilde{U}^-(\infty, w^*_j) = \overline{U^-(\infty,w^*_j) - (S^-(w^*_j) \cup
S^+(w^*_{j+1}))}$. Then branches of $\log(\pi(w) - \pi(w^*_j)),
\log(\pi(w))$ map $\tilde{U}(w^*_j), \tilde{U}^-(\infty, w^*_j)$
univalently to semi-infinite horizontal strips $V_j, W_j$ of the form $\{ \Re z \leq A_j, B_j \leq \Im
z \leq C_j \}$ and  $\{ \Re z \geq A'_j, B_j \leq \Im
z \leq C_j \}$ respectively. Similarly appropriate logarithms map $S^-(w^*_j), S^+(w^*_j)$ to
lower and upper half-planes $L_j, H_j$ of the form $\{ \Im z \leq
B_j \}$ and $\{ \Im z \geq C_j\}$ respectively.

\medskip

Clearly $W(\RR_i)$ is homeomorphic to the disjoint union of the
$V_j, W_j, L_j, H_j$ glued together as follows: each strip $V_j$ is
glued to the half-planes $L_j, H_j$ along the boundary arcs $\{
\Re z \leq A_j, \Im z = B_j \}$ and $\{ \Re z \leq A_j, \Im z = C_j \}$
respectively, while each strip $W_j$ is glued to the half-planes $L_j,
H_{j+1}$ along the boundary arcs $\{ \Re z \geq A'_j, \Im z = B_j \}$
and $\{ \Re z \geq A_j, \Im z = C_{j+1} \}$ respectively. It is
easy to see that the resulting quotient is homeomorphic to a
closed disk with a puncture $\{ 0 < |z| \leq 1 \}$. $\diamond$

\medskip

\begin{prop} \label{topology} The Riemann surface $\SSS^{\times}$
is homeomorphic to a closed surface with finitely many punctures
(one for each component of $E$).
\end{prop}

\medskip

\noindent{\bf Proof:} It follows from Lemmas
\ref{nbhdinfinity},\ref{ends1} and \ref{ends4} above that
$\overline{\SSS^{\times} - E}$ is a compact surface with finitely
many boundary components which are Jordan curves, and $\SSS^{\times}$ is
given by attaching a closed disk with a puncture to each boundary components,
hence $\SSS^{\times}$ is a closed surface with finitely many punctures. $\diamond$

\medskip

\section{Uniformization theorems}

\medskip

We have seen above that a finite type log-Riemann surface has
finitely many ends, each homeomorphic to a punctured disk (henceforth by an
'end' we will mean a component of a neighbourhood of infinity). We
shall show that each is indeed biholomorphic to a punctured disk,
by showing that each end can be isometrically embedded in a
log-Riemann surface whose finite completion is simply connected,
and then using the uniformization theorem of \cite{bipm2} to
conclude that the end is biholomorphic to a neighbourhood of
infinity in $\C$, hence is a punctured disk.

\medskip

It will turn out that to each end is associated an integer $K$, the index of the end,
which corresponds to the index of a holomorphic vector field with an exponential singularity
associated to the end. We define below a countable family of log-Riemann surfaces $\SSS(w_0, \dots, w_{n-1}, w, K)$
indexed by an integer $K$ such that $\SSS(w_0, \dots, w_{n-1}, w, K)$ has an end with index $K$;
any end of a finite type log-Riemann surface with index $K$ will then be embeddable
into $\SSS(w_0, \dots, w_{n-1}, w, K)$.

\medskip

\subsection{A family of finite type log-Riemann surfaces}

\medskip

Given $n$ points $w_0$, $\dots$, $w_{n-1} \in \C$ (not necessarily distinct), a point $w$ distinct
from $w_0, \dots, w_{n-1}$ and an integer $K \in \integers$ we define a finite type
log-Riemann surface $\SSS(w_0, \dots, w_{n-1}, w, K)$ as follows:

\medskip

Choose a point $z_0$ not lying on
any of the lines passing through $w_i, w_j, 0 \leq i < j \leq n-1$,
and a point $w$ not lying on any of the lines passing through $z_0$ and
$w_j, j = 0,\dots, n-1$. Let $l, l_0, \dots, l_{n-1} \subset \C$ be
the closed half-lines starting at the points $w,w_0,\dots,w_{n-1}$ in the directions
$(w - z_0), (w_0 - z_0), \dots, (w_{n-1} - z_0)$ respectively.
Consider the slit planes $C := \C - l, C_j := \C - l_j, C^*_j
:= \C - (l \cup l_j), j = 0, \dots, n-1$. Completing these slit planes
with respect to the path-metric induced from $\C$ gives metric
spaces $\overline{C} = C \sqcup (l^+ \sqcup l^-)$,
$\overline{C_j} = C_j \sqcup (l^-_j \sqcup l^+_j)$, $\overline{C^*_j} = C^*_j \sqcup (l^- \sqcup l^+)
\sqcup (l^-_j \sqcup l^+_j)$ where $l^{\pm}, l^{\pm}_j$ are isometric copies of
$l, l_j$ respectively (representing the 'top' and 'bottom' sides of the lines). We construct the
log-Riemann surface $\SSS(w_0, \dots, w_{n-1}, w, K)$ by pasting together copies
of $\overline{C}, \overline{C_j}, \overline{C^*_j}$ along the 'slits' $l^{\pm}, l^{\pm}_j$
as follows:

\medskip

For $j = 0,\dots, n-1$ we take a copy of $\overline{C^*_j}$ and a
family $(\overline{C_j}^{(k)})_{k \in \integers}$ of copies
$\overline{C_j}^{(k)}$ of $\overline{C_j}$. We treat the cases $K=0, K > 0, K < 0$
separately:

\medskip

\noindent{\bf (1) $K = 0 :$} For $j = 0,\dots, n-1$, we identify
isometrically the following lines: $l^-
\subset \overline{C^*_j}$ is identified with $l^+ \subset
\overline{C^*_{j+1}}$, $l^-_j \subset \overline{C_j}^{(0)}$ with $l^+_j \subset
\overline{C^*_j}$, $l^-_j \subset \overline{C^*_j}$ with $l^+_j \subset
\overline{C_j}^{(1)}$, and $l^-_j \subset \overline{C_j}^{(k)}$ with $l^+_j \subset
\overline{C_j}^{(k+1)}$ for $k \neq 0$.

\medskip

It is not hard to see that we obtain a log-Riemann surface $\SSS = \SSS(w_0, \dots, w_{n-1}, w, 0)$
with $n$ ramification points of infinite order projecting onto the points
$w_0, \dots, w_{n-1}$ and one of order
$n$ projecting onto $w$, such that $\SSS^{\times}_i$ is simply connected. The
log-Riemann surface $\SSS$ has a skeleton with stars $C^*_j$ forming a cycle $C$ of
length $n$, each being attached to two half-lines formed by the stars
$(C^{(k)}_j)_{k \leq 0}$ and $(C^{(k)}_j)_{k \geq 1}$.

\medskip

\noindent{\bf (2) $K > 0 :$} In this case we take $K$
copies $(\overline{C}^{l})_{1 \leq l \leq K}$
of $\overline{C}$ as well. For $j = 0,\dots, n-2$, we make the same
isometric identifications as above:
$l^- \subset \overline{C^*_j}$ is identified with $l^+ \subset
\overline{C^*_{j+1}}$, $l^-_j \subset \overline{C_j}^{(0)}$ with $l^+_j \subset
\overline{C^*_j}$, $l^-_j \subset \overline{C^*_j}$ with $l^+_j \subset
\overline{C_j}^{(1)}$, and $l^-_j \subset \overline{C_j}^{(k)}$ with $l^+_j \subset
\overline{C_j}^{(k+1)}$ for $k \neq 0$. We make the following further isometric
identifications: $l^- \subset \overline{C^*_{n-1}}$ is identified
with $l^+ \subset \overline{C}^{(1)}$, $l^- \subset
\overline{C}^{(K)}$ with $l^+ \subset \overline{C^*_0}$, and if $K > 1$ then
$l^- \subset \overline{C}^{(i)}$ is identified
with $l^+ \subset \overline{C}^{(i+1)}$ for $i = 1, \dots, K - 1$.

\medskip

In this case we obtain a log-Riemann surface $\SSS(w_0, \dots, w_{n-1}, w, K)$
with $n$ ramification points of infinite order projecting onto the points
$w_0, \dots, w_{n-1}$ and one of order
$n+K$ projecting onto $w$, such that $\SSS^{\times}(w_0, \dots, w_{n-1}, w, K)$ is simply connected. The
log-Riemann surface has a skeleton with the stars
$C^*_0, \dots, C^*_{n-1}, C^{(1)}, \dots, C^{(K)}$
forming a cycle of length $n + K$, with each $C^*_j$ being attached to
two half-lines formed by the stars
$(C^{(k)}_j)_{k \leq 0}$ and $(C^{(k)}_j)_{k \geq 1}$.

\medskip

From the uniformization theorem of \cite{bipm2} we have:

\medskip

\begin{prop} \label{genuszerowithramfn} The finite completion
$\SSS^{\times}(w_0, \dots, w_{n-1}, w, K)$ for $K \geq 0$ is
biholomorphic to $\C$. There is a polynomial $P$ of
degree $n$ and a uniformization
$F : \C \to \SSS^{\times}(w_0, \dots, w_{n-1}, w, K)$ satisfying
$(\pi \circ F)'(z) = z^{n+K} e^{P(z)}$.
\end{prop}

\medskip

We remark that the $n$ infinite ramification points of $\SSS(w_0,
\dots, w_{n-1}, w, K)$ all belong to one cycle. The corresponding
end $W(\RR_{\infty})$ is the only end of \hfill
$\SSS^{\times}(w_0, \dots, w_{n-1}, w, K)$. Under the uniformization $F$ above, the end $W(\RR_{\infty})$
is biholomorphic to a punctured disk neighbourhood of $\infty$ in
$\hat{\C}$.

\medskip

\noindent{\bf (3) $K < 0 :$ } In this case we replace the two copies
$\overline{C_j}^{(1)}, \overline{C_j}^{(1 + (- K))}$ of
$\overline{C_j}$ with two additional copies $\overline{C^*_0}^{(1)},
\overline{C^*_0}^{(2)}$ of $\overline{C^*_0}$ instead. For $j = 0,\dots, n-1$,
we make the following isometric identifications:
$l^- \subset \overline{C^*_j}$ is identified with $l^+ \subset
\overline{C^*_{j+1}}$, $l^-_j \subset \overline{C_j}^{(0)}$ with $l^+_j \subset
\overline{C^*_j}$. For $j = 1,\dots, n-1$ we identify
$l^-_j \subset \overline{C^*_j}$ with $l^+_j \subset \overline{C_j}^{(1)}$,
and $l^-_j \subset \overline{C_j}^{(k)}$ with
$l^+_j \subset \overline{C_j}^{(k+1)}$ for $k \neq 0$. We
identify $l^-_0 \subset \overline{C^*_0}$ with
$l^+_0 \subset \overline{C^*_0}^{(1)}$, $l^-_0 \subset \overline{C^*_0}^{(2)}$ with
$l^+_0 \subset \overline{C_0}^{(2 + (-K))}$, $l^+ \subset
\overline{C^*_0}^{(1)}$ with $l^- \subset \overline{C^*_0}^{(2)}$,
$l^- \subset \overline{C^*_0}^{(1)}$ with $l^+ \subset
\overline{C^*_0}^{(2)}$. If $K < -1$ then we identify $l^-_0
\subset \overline{C^*_0}^{(1)}$ with $l^+_0 \subset
\overline{C_0}^{(2)}$, $l^-_0 \subset \overline{C_0}^{(-K)}$ with
$l^+_0 \subset \overline{C^*_0}^{(2)}$ and if $K < -2$ then we
identify $l^-_0 \subset
\overline{C_0}^{(i)}$ with $l^+_0 \subset
\overline{C_0}^{(i+1)}$ for $i = 2, \dots, (-K)-1$. Finally we
identify $l^-_0 \subset \overline{C_0}^{(k)}$ with
$l^+_0 \subset \overline{C_0}^{(k+1)}$ for $k \leq -1$ and $k \geq
2+(-K)$.

\medskip

In this case the log-Riemann surface $\SSS(w_0, \dots, w_{n-1}, w, K)$
we obtain is of finite type, has again $n$ infinite order ramification points
projecting onto $w_0, \dots, w_{n-1}$, but now two finite order
ramification points projecting onto $w$, one of order $n$ and one
of order two. Moreover there is an end $U$ of the finite completion
which is a finite-sheeted covering of $\{ |z| > R \}$ of degree
$(-K)$. The finite completion is {\it not} simply connected in
this case, but we still have:

\medskip

\begin{prop} \label{genuszerowithpole} The finite completion
$\SSS^{\times}(w_0, \dots, w_{n-1}, w, K)$ for $K < 0$ is
biholomorphic to the punctured plane $\C^*$. There is a uniformization
$F : \C^* \to \SSS^{\times}(w_0, \dots, w_{n-1}, w, K)$ satisfying
$(\pi \circ F)'(z) = z^K (z - p_1)^n (z - p_2) e^{P(z)}$ for some polynomial $P$ of
degree $n$ and two distinct non-zero $p_1,p_2 \in \C$.
\end{prop}

\medskip

\noindent{\bf Proof:} It is straightforward to see that the
finitely completed skeleton $\Gamma^{\times}$ of $\SSS = \SSS(w_0, \dots,
w_{n-1}, K)$ has only one cycle corresponding to the end $U$ which is
a degree $(-K)$ covering of $\{|z| > R \}$. The end $U$ is
biholomorphic to a punctured disk, hence we may add a point $q$ to
$\SSS^{\times}$ to obtain a Riemann surface
$\SSS^*$ such that $U \cup \{q\}$ is biholomorphic to a disk, and
$\pi$ extends to have a pole of order $(-K)$ at $q$. Moreover the
surface $\SSS^*$ is simply connected. Since it has only finitely many
ramification points, using the Kobayashi-Nevanlinna parabolicity
criterion of \cite{bipm2} it follows that $\SSS^*$ is biholomorphic to $\C$.

\medskip

We now argue as in the proof of Theorem 1.1 of \cite{bipm2}. Let $F : \C \to \SSS^*$
be a uniformization such that $F(0) = q$, and choose a basepoint $q' \neq q$ in $\SSS$.
The approximation Theorem 2.10 of \cite{bipm2} gives a sequence of
pointed finite-sheeted log-Riemann surfaces $(\SSS_k, q'_k)$
converging to $(\SSS, q')$ in the sense of Caratheodory (see
\cite{bipm1}). For $k$ large the surfaces
$\SSS^{\times}_k$ also have one end which is a degree $(-K)$ covering
of $\{|z| > R \}$, and hence as above we can add a point $q_k$ to
obtain a Riemann surface $\SSS^*_k$ which is biholomorphic to
$\C$. Let $F_k : \C \to \SSS^*_k$ be a uniformization such that
$F_k(0) = q_k$. Since $\pi_k \circ F_k$ is finite-to-one in a neighbourhood
of infinity, it follows that $\pi_k \circ F_k$ is a rational function, which
moreover has only one finite pole of order $(-K)$ at $0$, two critical
points of orders $n-1$ and $1$ corresponding to the finite ramification points of $\SSS_k$ projecting
onto $w$, and $n$ critical points of increasing orders (as $k \to \infty$) corresponding
to the finite ramification points projecting onto $w_0, \dots, w_{n-1}$.
By the Caratheodory convergence theorem of \cite{bipm2},
normalizing the $F_k$'s appropriately, we have $\pi_k \circ F_k \to \pi \circ F$
uniformly on compacts of $\C^*$ as $k \to \infty$.

\medskip

Now the same argument as in \cite{bipm2} shows that the
nonlinearity of $\pi \circ F$ is a rational
function of the form
$$
\frac{(\pi \circ F)''(z)}{(\pi \circ F)'(z)} = \frac{K-1}{z} + \frac{n-1}{z - p_1} + \frac{1}{z - p_2} + P'(z)
$$
for two distinct non-zero $p_1, p_2 \in \C$ and some polynomial $P$ of degree $n$, and the
result follows upon integration of the above. $\diamond$

\medskip

We note that as before the $n$ infinite ramification points
all lie in the same cycle. However in this case ($K < 0$), the surface
$\SSS^{\times}(w_0, \dots, w_{n-1}, w, K)$ has two ends, one of the form
$W(\RR_{\infty})$ corresponding to the cycle of infinite ramification points,
and the other the degree $(-K)$ covering $U$ of $\{|z| \geq R \}$.
Under the uniformization $F : \C^* \to \SSS^{\times}(w_0, \dots,
w_{n-1}, w, K)$ these are biholomorphic to punctured disks in $\C$.

\medskip

\subsection{Isometric embedding of ends}

\medskip

Let $\SSS$ be a log-Riemann surface of finite type and
consider the decomposition of the skeleton
$$
\Gamma = \Gamma_0 \sqcup \left(\sqcup_{w^* \in \RR_{\infty}} L^+(w^*) \sqcup
L^-(w^-)\right).
$$
We observe that adding to $\Gamma_0$ finite initial
segments of each half-line $L^{\pm}(w^*), w^* \in \RR_{\infty}$ gives
a similar decomposition for which all properties proved in the previous
section still hold. Choosing these initial segments to be large enough
and of equal lengths,
we may arrange that the decomposition satisfies the
following condition:

\medskip

The boundary curves $\partial U(w^*), \partial U^-(\infty, w^*),
w^* \in \RR_{\infty}$ define bi-infinite lines in the skeleton,
$\Gamma(w^*)$ (defined earlier) and $\Gamma'(w^*)$ (say)
respectively. We have $(L^-(w^*) \cup L^+(w^*)) \subset \Gamma(w^*),
(L^-(w^*) \cup L^+(u(w^*))) \subset \Gamma'(w^*)$,
and $\Gamma(w^*) - (L^-(w^*) \cup L^+(w^*)), \Gamma'(w^*) - (L^-(w^*) \cup L^+(u(w^*)))$
are finite segments. We will assume that $\Gamma_0$ has been
chosen so that all these segments have lengths lying in
an interval $[2N - c_1, 2N + c_1]$ where $N, c_1, c_2 \geq 2$ are integers such that
$N \geq 8(\#\RR_{\infty}+1)(c_1 + c_2)$.

\medskip

Now let $\RR_i = \{ w^*_0, \dots, w^*_{n-1} \} \subset \RR_{\infty}$
be a cycle of infinite ramification points (with $w^*_{j+1} =
u(w^*_j)$), and $W(\RR_i)$ the corresponding component of $E$. Let $w_j = \pi(w^*_j)$.
We have:

\medskip

\begin{theorem} \label{embedding} There exists $K \in \integers$ and an
isometric embedding $\iota : W(\RR_i) \to \SSS(w_0, \dots, w_{n-1}, w, K)$
such that $\tilde{\pi} \circ \iota = \pi$ and $\iota(W(\RR_i))$ is
the end of $\SSS^{\times}(w_0, \dots, w_{n-1}, w, K)$ corresponding to
the cycle of infinite ramification points of $\SSS^{\times}(w_0, \dots, w_{n-1}, w, K)$
(here $\tilde{\pi}$ is the projection mapping of $\SSS(w_0, \dots, w_{n-1}, w, K)$).
\end{theorem}

\medskip

\noindent{\bf Proof:} We first show that (in the notation of Lemma
\ref{ends4}) the union of all but one of the sets
making up $W(\RR_i)$, namely $S^-(w^*_0)$, $\tilde{U}^-(\infty, w^*_0)$, $S^+(w^*_1)$,
$\tilde{U}(w^*_1)$, $S^-(w^*_1), \dots, S^-(w^*_{n-1})$,
$\tilde{U}^-(\infty, w^*_{n-1})$,$S^+(w^*_0)$, can be
isometrically embedded into $\SSS(w_0, \dots, w_{n-1}, K)$ for $|K| \leq N$,
and the embedding then extends to the remaining 'piece' $\tilde{U}(w^*_0)$ if
$K$ is chosen appropriately.

\medskip

Let $|K| \leq N/2$ and let $\SSS_K = \SSS(w_0, \dots, w_{n-1}, w,
K)$. We recall that $\SSS_K$ contains a family
$(\overline{C_0}^{(-k)})_{k \geq 1}$ of copies of
$\overline{C_0}$. Given $k_0 \geq 1$, clearly there is a unique isometry $\iota$ of
$S^-(w^*_0)$ onto $\cup_{k \geq k_0} \overline{C_0}^{(-k)}$ such that $\pi_i \circ \iota =
\pi$. We first show that $k_0$ can be chosen so
that $\iota$ extends to an isometric embedding of the union of
$S^-(w^*_0)$, $\tilde{U}^-(\infty, w^*_0)$, $S^+(w^*_1)$,
$\tilde{U}(w^*_1)$, $S^-(w^*_1), \dots, S^-(w^*_{n-1})$,
$\tilde{U}^-(\infty, w^*_{n-1})$, and $S^+(w^*_0)$. The indices $j, j+1$
will be taken modulo $n$ throughout.

\medskip

Suppose for some $0 \leq j \leq n-2$ we are given $k_j \geq 1$ and
an isometric embedding $\iota$ of $S^-(w^*_j)$ into $\SSS_K$ with
image $\cup_{k \geq k_j} \overline{C_j}^{(-k)}$.
If $1 \leq k_j \leq 2N - c_1 - 2$ then there is a unique isometric extension of $\iota$ to
$\tilde{U}^-(\infty, w^*_j)$, with image passing through the stars
$C^{(-k_j+1)}_j, \dots, C^{(0)}_j$, $C^*_j$, $C^{(1)}_{j+1}, \dots,
C^{(k'_{j+1})}_{j+1}$, where $k'_{j+1} = a'_j - (k_j + 1) \geq 1$, $a'_j$ being the
length of the segment $\Gamma'(w^*_j) - (L^-(w^*_j) \cup L^+(w^*_{j+1})))$
(in the skeleton $\Gamma$). This determines a further isometric
extension to $S^+(w^*_{j+1})$ with image in the stars $C^{(k)}_{j+1},
k \geq k'_{j+1} + 1$.

\medskip

There is also an isometric
extension to $\tilde{U}(w^*_{j+1})$ which, if $k'_{j+1} \leq 2N - c_1 - 1$,
has image in the stars $C^{(k'_{j+1})}_{j+1}, \dots, C^{(1)}_{j+1}, C^*_{j+1},
C^{(0)}_{j+1}, \dots, C^{(-k_{j+1})}_{j+1}$ where $k_{j+1} = a_{j+1} - k'_{j+1} \geq 1$,
$a_{j+1}$ being the length of the segment
$\Gamma(w^*_{j+1}) - (L^-(w^*_{j+1}) \cup L^+(w^*_{j+1}))$. Again, this gives
a unique isometric extension to $S^-(w^*_{j+1})$ with image in the
stars $C^{(-k)}_{j+1}, k \geq k_{j+1} + 1$.

\medskip

By hypothesis $a'_j \in [2N-c_1, 2N + c_1]$ so if $k_j \in [N - jc_1 - c_2, N +
jc_1 + c_2]$ then
$k'_{j+1} = a'_j - (k_j + 1) \in [N - (j+1)c_1 - c_2 - 1, N + (j+1)c_1 + c_2 + 1] \subset [1, 2N - c_1 - 2]$,
and $k_{j+1} = a_{j+1} - k'_{j+1} \in [N - (j+1)c_1 - c_2, N + (j+1)c_1 + c_2] \subset [1, \infty)$
(using the hypotheses on $N, c_1, c_2$).
It follows that choosing $k_0 \in [N - c_2, N + c_2]$
we can extend $\iota$ inductively to an isometry defined
on the union of $S^-(w^*_0)$, $\tilde{U}^-(\infty, w^*_0)$, $S^+(w^*_1)$,
$\tilde{U}(w^*_1)$,$S^-(w^*_1)$, $\dots$, $S^-(w^*_{n-1})$.

\medskip

The image of $S^-(w^*_{n-1})$ lies in the stars $C^{(-k)}_{n-1}, k
\geq k_{n-1}+1$, and $k_{n-1} \leq N + (n-1)c_1 + c_2$. As
$\tilde{U}^-(\infty, w^*_{n-1})$ corresponds to a segment in the
skeleton of length $a'_{n-1} \geq 2N - c_1$ and $|K| \leq N/2$,
it is not hard to see that $\iota$ has a unique isometric extension to
$\tilde{U}^-(\infty, w^*_{n-1})$ with image starting in the star
$C^{(-k_{n-1})}_{n-1}$ and ending in a star $C^{(k'_0)}_0$, where
$k'_0 \geq 1$ is given by $k'_0 = a'_{n-1} - (k_{n-1}+1+K)$. Then
$\iota$ has a unique isometric extension to
$S^+(w^*_0)$, with image lying in the stars $C^{(k)}, k \geq k'_0 +
1$.

\medskip

It remains to define $\iota$ on $\tilde{U}(w^*_0)$.
As $\iota$ is given on $\tilde{U}(w^*_0) \cap S^+(w^*_0)$ and $S^-(w^*_0)$,
there is an extension to $\tilde{U}(w^*_0)$ if and only if the
number of stars from $C^{(k'_0)}_0$ to $C^{(-k_0)}$ is equal to the length
of the segment $\Gamma(w^*_0) - (L^-(w^*_0) \cup L^+(w^*_0))$, or in other words
if $k'_0 + k_0 + 1 = a_0$.

\medskip

Using the recursion formulae for $k_j, k'_j$, this reduces to
$$
K = \sum_{j = 0}^{n-1} (a'_j - a_j) - (n-1)
$$
and by the hypothesis on the lengths $a'_j, a_j$ the right-hand
side above is bounded above by $2nc_1+(n-1) \leq N/2$. Hence $K$
can be chosen as required so that $\iota$ extends to an isometric
embedding of $W(\RR_i)$ into $\SSS_K$. $\diamond$

\medskip

\begin{cor} \label{endsdisks} Each end of a finite type
log-Riemann surface is biholomorphic to a punctured disk.
\end{cor}

\medskip

\noindent{\bf Proof:} Each end is either a component $U$ of
$U(\infty)$ such that $\pi : U \to \{ |z| \geq R \}$ is a
finite-sheeted covering, in which case the interior of $U$ is
biholomorphic to a punctured disk, or is a component of $E$ of the
form $W(\RR_i)$. It follows from the previous Theorem and the remarks
following Propositions \ref{genuszerowithramfn},
\ref{genuszerowithpole} that the end $W(\RR_i)$ is biholomorphic
to a punctured disk.
$\diamond$

\medskip

\subsection{Proofs of main theorems}

\medskip

We are now in a position to prove the main Theorems
\ref{uniformthm}, \ref{conversethm}:

\medskip

\noindent{\bf Proof of Theorem \ref{uniformthm}:}  It follows immediately from
Proposition \ref{topology} and Corollary \ref{endsdisks} that the finite completion
$\SSS^{\times}$ of a finite type log-Riemann surface is
biholomorphic to a closed Riemann surface with finitely many punctures
corresponding to the ends of $\SSS^{\times}$.
For an end which is a finite-sheeted degree $d$ covering of a
neighbourhood of $\infty$ in $\hat{\C}$, it follows that
$\pi$ extends meromorphically to have a pole of order $d$ at the corresponding puncture.

\medskip

For an end $W(\RR_i)$ corresponding to a cycle of infinite
ramification points, by Theorem \ref{embedding} and Propositions
\ref{genuszerowithramfn}, \ref{genuszerowithpole}, we obtain a
biholomorphic map $F$ from a punctured disk neighbourhood of $\infty$
in $\hat{\C}$ to $W(\RR_i)$ such that in terms of the local coordinate $z = F^{-1}$,
$d\pi = R(z) e^{P(z)} dz$ for some polynomial $P$ and some rational function
$R$. It follows that $d\pi$ has an exponential singularity at the
corresponding puncture. $\diamond$

 $\diamond$

\medskip

\noindent{\bf Proof of Theorem \ref{conversethm}:} Given a non-constant
meromorphic map $f$ on a closed Riemann surface with punctures $S$ such that
$df$ has exponential singularities at the punctures, we note that $df$
has no poles or zeroes in a neighbourhood of each puncture, hence the set $A$
of poles and critical points of $f$ is finite. Let $S'$ be the punctured surface
$S' = S - A$. Then $f : S' \to \C$ is a local diffeomorphism, and it suffices
to show that finitely many points are added in the completion $\overline{S}$
of $S'$ with respect to the path metric induced by $|df|$. If $p \in S'$ tends
to a puncture which is a critical point of $f$, then $p$ tends to a unique limit in
$\overline{S}$, which is a finite ramification point of $S'$, while if
$p$ tends to a pole of $f$ then $p$ tends to infinity in $\overline{S}$.

\medskip

It remains to show that as $p$ tends to a puncture $p_0$ which is an
exponential singularity of $df$, then $p$ can only accumulate in
$\overline{S}$ on a finite set of points (only depending on $p_0$).
Let $df = e^h \omega$ near $p_0$, where $h$ is a meromorphic
function with a pole of order $n$ say at $p_0$ and $\omega$ is a $1$-form
meromorphic near $p_0$. We can choose a local coordinate $z$ such
that $z(p) = \infty, h(z) = z^n$, and $df = g(z) e^{z^n} dz$
where $g$ is a function meromorphic near $z = \infty$. Then for
some integer $k$ and some $C > 0$ we have
$$
\frac{1}{C} |z^k e^{z^n}| |dz| \leq |df| \leq C |z^k e^{z^n}| |dz|
$$
It follows that the path-metrics $d, d'$ induced by $|df|$ and $|z^k
e^{z^n}| |dz|$ on a punctured neighbourhood $D$ of $p_0$
are bi-Lipschitz equivalent, hence so are the completions of
$D$ with respect to $d, d'$.

\medskip

If $k \geq 0$ then by Theorem 1.2 of \cite{bipm2}
the function $\int z^k e^{z^n} dz$ defines a log-Riemann surface with $n$ infinite
ramification points and we are done.

\medskip

If $k  = -m < 0$ then the metric induced by $|z^k e^{z^n}dz|$ is bi-Lipschitz
equivalent to that induced by $|\eta| = |(1/z^m - C/z^{m+n})e^{z^n} dz|$ where
the constant $C$ is chosen so that the residue of the $1$-form $\eta$ at $\infty$
vanishes, so $\eta$ has a primitive $F = \int \eta$ on $\C^*$. We can approximate
$F$ by rational functions
$$
R_N(z) = \int \left(\frac{1}{z^m} -
\frac{C_N}{z^{m+n}}\right)\left(1+ \frac{z^n}{N}\right)^N dz
$$
where the constants $C_N \to C$ are chosen so that the residue at $\infty$ of
the $1$-form in the integral vanishes. Then each $R_N$ defines a log-Riemann
surface structure $\SSS_N$ on $\C^*$
with $n$ finite order ramification points each of order
$N+1$. By the compactness Theorem 2.11 of \cite{bipm2}, there is a subsequence
of $\SSS_N$ which converges in the sense of Caratheodory to a log-Riemann
surface structure $\SSS$ on $\C^*$
with at most $n$ ramification points. By the Caratheodory convergence
Theorem 1.2 of \cite{bipm1}, the maps $R_N$ converge along the same subsequence
uniformly on compacts to the projection mapping of $\SSS$, which is
hence given by $F$. It follows that there are at most $n$ points
added in the completion of a neighbourhood of $\infty$ with respect to
the metric induced by $|\eta|$. $\diamond$

$\diamond$

\medskip

\bibliography{unifhighgenus}
\bibliographystyle{alpha}

\end{document}